\documentclass[reqno,11pt]{amsart}
\usepackage{latexsym}
\usepackage{fullpage}
\usepackage{amssymb}
\def\bea{\begin{eqnarray}}
\def\eea{\end{eqnarray}}
\def\bdm{\begin{displaymath}}
\def\edm{\end{displaymath}}
\def\R{\ensuremath{\mathbb R}}
\def\C{\ensuremath{\mathbb C}}

\def\Z{\ensuremath{\mathbb Z}}

\def\S{\mathcal S}
\def\k{\kappa}

\def\id{\mbox{id}}

\def\rk{\operatorname{rk}}
\def\ind{\operatorname{ind}}
\def\Spin{\operatorname{Spin}}
\def\tr{\operatorname{tr}}
\def\Adach{{\hat A}}
\def\rk{\operatorname{rk}}

\def\Sign{\operatorname{Sign}}
\def\Sym{\operatorname{Sym}}
\def\ch{\operatorname{ch}}

\def\tensor{\otimes}
\def\even{{\operatorname{even}}}
\def\odd{{\operatorname{odd}}}
\def\Cl{\mathop{\mathrm C\mskip-1.5mu\ell\mskip1.5mu}}

\numberwithin{equation}{section}
\swapnumbers
\theoremstyle{plain}
\newtheorem{Lemma}[equation]{Lemma}
\newtheorem{Proposition}[equation]{Proposition}
\newtheorem{Corollary}[equation]{Corollary}
\newtheorem{Theorem}[equation]{Theorem}
\newtheorem{Question}[equation]{Question}

\theoremstyle{definition}

\theoremstyle{remark}
\newtheorem{Remark}[equation]{Remark}

\begin{document}

\title{Scalar Curvature Estimates for Compact Symmetric Spaces}
\author{S. Goette, U. Semmelmann}
\address{Math.\ Inst.\ Uni.\ T\"ubingen, Auf der Morgenstelle~10,
D-72076 T\"ubingen, Germany}
\email{sebastian.goette@uni-tuebingen.de}
\address{Math.\ Inst.\ Uni.\ M\"unchen, Theresienstr.~39, D-80333 M\"unchen,
Germany}
\email{semmelma@rz.mathematik.uni-muenchen.de}
\thanks{Both authors were supported by a research fellowship of the DFG}
\begin{abstract}
We establish extremality of Riemannian metrics~$g$ with non-negative 
curvature operator on symmetric spaces~$M=G/K$ of compact type 
with~$\rk G-\rk K\le 1$. Let~$\bar g$ be another metric with scalar 
curvature~$\bar\kappa$, such that~$\bar g\ge g$ on 2-vectors.
We show that~$\bar\kappa\ge\kappa$ everywhere on~$M$
implies~$\bar\kappa=\kappa$. Under an additional condition on the 
Ricci curvature of~$g$, $\bar\kappa\ge\kappa$ even implies~$\bar g=g$.
We also study area-non-increasing spin maps onto such 
Riemannian manifolds.
\end{abstract}

\maketitle

There is a well known relation between the existence of
metrics of positive scalar curvature on a compact manifold
$M$ and the topology of $M$. Given such a metric $g$ with scalar
curvature $\kappa > 0$, it is interesting to ask how large~$\kappa$ 
can become as a function on~$M$ when one varies the metric~$g$. 
Of course one should not
allow scaling of the metric, so one has to compare $g$ with
suitable other metrics, e.g. with metrics $\bar g$ which do not decrease
areas with respect to the fixed metric $g$. These are the
metrics in
\begin{equation}\label{MgDef}
	{\mathcal M}(g)
	:=\bigl\{\,{\bar g}\in \Sym^2(TM)\bigm|
		\left|v\wedge w\right|_{\bar g}\ge\left|v\wedge w\right|_g
		\text{for all~$v$, $w \in TM$}\,\bigr\}\;.
\end{equation}
For ${\bar g} \in {\mathcal M}(g)$ we will write ${\bar g} \ge g$ on 2-vectors.
Using the K-area inequalities, M.~Gromov showed 
in \cite{gromov} that there is a finite upper bound for the
minimum of the scalar curvature if one varies over metrics
in ${\mathcal M}(g)$. However it remains a problem to find
sharp upper bounds in terms of the curvature of the fixed
metric $g$. A first example for a sharp upper bound was given
by M.~Llarull in \cite{ll2}. He considered metrics on the 
sphere $(S^n,g)$ where $g$ is the metric of constant curvature.
If $\bar g$ is any metric on $S^n$ with scalar curvature $\bar \kappa$
and  with ${\bar g} \ge g$ on 2-vectors, he showed that 
(i)
$\bar\kappa(p)\le\kappa(p)$ for some~$p\in M$
and 
(ii) ${\bar \kappa}\ge \kappa$ implies ${\bar \kappa}=
\kappa$.
Indeed he showed that~(iii) 
$\bar\kappa \ge \kappa$ even implies~${\bar g} = g$.
 
Metrics $g$ having property (ii) will be called {\it area-extremal}.
Note that~(i) follows from~(ii).
If $g$ is an area-extremal metric of constant scalar curvature, this constant provides a sharp upper bound for the minimum of the scalar curvature of all metrics in 
${\mathcal M}(g)$. 

Let us relate area-extremality to a theorem of J.~Lohkamp,
which sharpens earlier results of J.~L.~Kazdan and F.~W.~Warner (\cite{KW}).
Let~$g$ be a Riemannian metric on a manifold~$M$ with~$\dim(M)\ge\penalty1000\relax3$,
and let~$\kappa_0\colon M\to\R$ be any function such that~$\kappa_0\le\kappa$
everywhere on~$M$.
Then by~\cite{Lohkamp},
there exists a metric~$\bar g$,
which is $C^0$-close to~$g$,
such that~$\bar\kappa$ is $C^0$-close to~$\kappa_0$.
In particular, it is always possible to decrease both the metric and the
scalar curvature simultaneously.
On the other hand, if $g$ is area-extremal,
then by condition~(ii) one cannot simultaneously increase both $g$
and $\kappa$.
We apply the construction of J.~Lohkamp to see that not all metrics~$g$ on~$M$
are area-extremal if~$\dim M\ge 3$.
Let us start with an arbitrary metric~$\bar g$ on~$M$.
Using~\cite{Lohkamp}, we can construct a metric~$g$
which is $C^0$-close to~$\frac{1}{2}\,\bar g$
and has~$\kappa\ll\bar\kappa$.
Since clearly~$\bar g\in{\mathcal M}(g)$,
$g$ is not area-extremal.
If~$M$ admits a metric of positive scalar curvature,
we start with~$\bar\kappa>0$.
Then~$g$ can be chosen such that~$\kappa$
approximates~$\frac{1}{2}\,\bar\kappa$\
in the $C^0$-topology.
This shows that there are metrics with positive scalar curvature
that are not area-extremal.

In \cite{gromov}, M.~Gromov asked which manifolds possess area-extremal 
metrics and how such metrics may look like. He conjectured that Riemannian 
symmetric spaces should have area-extremal metrics. He also proposed to investigate not only variations of the metric on~$M$ itself, but to consider also area-non-increasing spin maps of non-vanishing $\Adach$-degree from other Riemannian manifolds to~$M$.

As mentioned above, M.~Llarull showed that the standard metric on $S^n$
is area-extremal with the additional rigidity (iii). 
In \cite{ll1},
he shows that~$\bar\kappa\ge\kappa\circ f$ implies~$\bar\kappa=\kappa\circ f$
if~$f$ is an area-non-increasing spin map of non-vanishing $\hat A$-degree
from a Riemannian manifold~$(N,\bar g)$ onto the round sphere. 
Later, M.~Min-Oo proved  that Hermitian symmetric spaces of compact type 
are area-extremal (cf. \cite{MinOo}). Finally, W.~Kramer proved in 
\cite{Wolfram} that quaternionic projective spaces are length-extremal,
which is a slightly weaker notion than area-extremality. In particular,
area-extremality implies length-extremality. 

In this paper,
we generalize the preceding results.
We prove area-extremality and rigidity for a certain class
of Riemannian metrics with non-negative curvature operator
on~$\Lambda^2(TM)$,
see Theorem~\ref{MainTh} below.
For this, we also require that either the Euler characteristic~$\chi(M)$ is non-zero,
or that a certain mod~$2$-index,
which is related to the Kervaire semi-characteristic,
does not vanish.

By the following theorem,
a compact,
simply connected Riemannian manifold with non-negative curvature operator
is homeomorphic to a symmetric space:

\begin{Theorem}[\cite{gal}, \cite{CC}, \cite{T}]\label{GMCTheorem}
If $(M,g)$ is a compact irreducible Riemannian manifold with 
non-negative curvature operator, then one of the following
cases must occur:
\begin{enumerate}
\item
the universal covering of $M$ is homeomorphic to a sphere, 
\item
the universal covering of $M$ is K\"ahler and biholomorphic
to a complex projective space,
\item
$M$ is locally symmetric.
\end{enumerate}
\end{Theorem}

The Euler characteristic of a Riemannian symmetric space $G/K$ of compact type
is different from zero iff $\rk(G)=\rk(K)$.
If the mod~2-index mentioned above is non-zero,
then~$\rk G-\rk K\le1$ and~$\dim M\equiv 0$, $1$ mod~$4$.
A certain stabilization trick allows us to treat all compact symmetric
spaces with~$\rk G-\rk K\le1$.
Let us summarize our main results,
which are obtained in Theorems~\ref{EulerSymmTh}
and~\ref{OddTh}.

\begin{Theorem}\label{MainTh}
Let~$(M,g)$ be a compact, connected, oriented Riemannian manifold
with non-negative curvature operator 
on~$\Lambda^2(TM)$,
such that the universal covering of~$M$ is homeomorphic
to a symmetric space~$G/K$ of compact type with~$\rk G\le\rk K+1$.
Let~$\bar g\in\mathcal M(g)$,
then~$\bar\kappa\ge\kappa$ implies~$\bar\kappa=\kappa$.
If moreover, the Ricci curvature of~$g$ satisfies~$\rho>0$
and~$2\rho-\kappa<0$,
then~$\bar\kappa\ge\kappa$ implies~$\bar g=g$.
\end{Theorem}

The conditions on~$\rho$ will be motivated at the end of this preface,
where we also conjecture a generalization of Theorem~\ref{MainTh}.

If~$\rk G=\rk K$,
we also compare~$\kappa$ with the scalar curvatures of metrics
on a different Riemannian manifold~$N$
via spin maps of non-vanishing $\hat A$-degree.
In Theorem \ref{MapTh} we prove again extremality and rigidity.

The proofs are based on a combination of the Bochner-Lichnerowicz-Weitzenb\"ock (BLW)
formula with the Atiyah-Singer index theorem applied to certain twisted Dirac operators.
We also need an estimate for the curvature term in the BLW formula
that uses non-negativity of the curvature operator.
In the odd-dimensional case ($\rk G=\rk K+1$)
we use the decomposition of the spinor bundle of~$G/H$.

In~\cite{GS},
we have established a similar result for K\"ahler manifolds
of positive Ricci curvature using different estimates.

The rest of the paper is organized as follows: 
In Section \ref{EstSect}, 
we investigate the BLW formula for a certain twisted Dirac operator.
The main result is contained in 
Lemma~\ref{Estimate}.
In Section \ref{InSect}, we apply 
this result. We use the index theorem in various settings to show the existence
of harmonic spinors.

We would like to thank  Ch.~B\"ar and G.~Weingart for
helpful comments and continued interest in our work.

\medskip\noindent
{\em Locally Area-Extremal Metrics.\/}
In the rest of this paper,
we generally consider globally area-extremal metrics on~$M$.
Here, ``globally'' means that~$\bar\kappa\ge\kappa$
implies~$\bar\kappa=\kappa$
for all metrics~$\bar g\in\mathcal M(g)$.
Here, we want to give a sufficient condition for a metric~$g$
to be locally area-extremal in the following sense:
There exists a neighborhood~$U(g)$ in the space~$\mathcal G$
of all Riemannian metrics on~$M$,
equipped with the $C^2$-topology,
such that~$\bar\kappa\ge\kappa$ implies~$\bar\kappa=\kappa$
for all metrics~$\bar g\in U(g)\cap\mathcal M(g)$.

\begin{Lemma}\label{LocalLemma}
Let~$g$ be a metric on a compact Riemannian manifold~$M$
whose Ricci curvature~$\rho$ is positive definite.
Then there exists a neighborhood~$U(g)$ of~$g$ in~$\mathcal G$
such that~$\bar\kappa\ge\kappa$ implies~$\bar\kappa=\kappa$
for all metrics~$\bar g\in U(g)$
with~$\left|v\right|_{\bar g}\ge\left|v\right|_g$ for all~$v\in TM$.

Suppose moreover that~$2\rho-\kappa$ is negative definite.
Then there exists another neighborhood~$U'(g)\subset\mathcal G$ of~$g$
such that~$\bar\kappa\ge\kappa$ implies~$\bar\kappa=\kappa$
for all metrics~$\bar g\in U'(g)\cap\mathcal M(g)$.
\end{Lemma}

\begin{proof}
Let~$\bar g$ be another metric on~$M$,
then there exists a $g$-symmetric endomorphism~$A$ of~$TM$,
such that
	$$\bar g(\mathord\cdot,\mathord\cdot)
	=g\bigl(e^A\mathord\cdot,\mathord\cdot\bigr)\;.$$
We consider the family
	$$g_t(\mathord\cdot,\mathord\cdot)
	=g\bigl(e^{tA}\mathord\cdot,\mathord\cdot\bigr)$$
with scalar curvature~$\kappa_t$ and Ricci curvature~$\rho_t$.
By a straightforward calculation in normal coordinates
around a point~$p$ in~$M$,
one checks that the derivative of~$\kappa_t$ is given by
\begin{equation}\label{KappaAbl}
	\frac\partial{\partial t}\kappa_t
	=\Delta_t(\tr_{g_t}A)
		-g_t\bigl((\nabla^{t,2}_{e_ie_j}A)e_i,e_j\bigr)
		-\tr_{g_t}\rho_t(A\mathord\cdot,\mathord\cdot)\;.
\end{equation}
Here, $\Delta_t$ denotes the Laplacian with respect to~$g_t$,
$\nabla^{t,2}$ denotes the second covariant derivative,
and~$\tr_{g_t}$ denotes the trace of a two-form with respect to~$g_t$.
Note that the first two terms can be written as divergences
with respect to the metric~$g_t$.
In particular, these terms either vanish identically,
or they become negative somewhere on~$M$.

Let us assume that the Ricci curvature~$\rho$ is positive definite.
If the metric~$\bar g$ is larger or equal than~$g$ on vectors,
then~$A$ is positive semi-definite.
In this case, the third term in~\eqref{KappaAbl} is~$\le 0$ everywhere on~$M$,
with equality iff~$A$ vanishes identically.
Because~$\rho$ depends continuously on the $2$-jet of~$g_t$,
there is some neighborhood~$U_1(g)$ in~$\mathcal G$
such that~$\rho'$ is positive definite for all~$g'\in U_1(g)$.
We define~$U(g)$ to be the set
of all metrics~$\bar g(\mathord\cdot,\mathord\cdot)
=g\bigl(e^A\mathord\cdot,\mathord\cdot\bigr)$
such that~$g\bigl(e^{tA}\mathord\cdot,\mathord\cdot\bigr)\in U_1(g)$
for all~$t\in[0,1]$.

If we have~$\bar g\ge g$ only on 2-vectors,
then the sum of any two eigenvalues of~$A$ is~$\ge 0$.
In other words, at most one eigenvalues~$a_i$ of~$A$ can be negative,
and its absolute value is not larger than any other eigenvalue.
On the other hand,
$2\rho-\kappa<0$ implies that no eigenvalue of~$\rho$ can be larger or equal
than the sum of the other eigenvalues.
In particular,
the condition~$2\rho_t-\kappa_t<0$ guarantees
that the last term in~\eqref{KappaAbl} is again non-positive,
with equality iff~$A$ vanishes identically.
Now, 
we define a neighborhood~$U_2(g)$ of~$g$ in~$\mathcal G$
such that~$\rho'>0$ and~$2\rho'-\kappa'<0$
for all~$g'\in U_1(g)$,
and we let~$U'(g)$ be the set
of all metrics~$g\bigl(e^A\mathord\cdot,\mathord\cdot\bigr)$
with~$A$ symmetric,
such that~$g\bigl(e^{tA}\mathord\cdot,\mathord\cdot\bigr)\in U_2(g)$
for all~$t\in[0,1]$.
\end{proof}

Note that in Lemma~\ref{LocalLemma},
we need precisely the same conditions on~$\rho$ as in Theorem~\ref{MainTh}.
One might even dare to ask the following

\begin{Question}\label{RicciQuestion}
Are all Riemannian metrics~$g$ on compact manifolds~$M$
with~$\rho>0$ and~$2\rho-\kappa<0$ area-extremal?
\end{Question}

For K\"ahler metrics, this has been answered affirmatively in~\cite{GS}.


\section{Scalar Curvature Estimates}\label{EstSect}
In this section,
we use the BLW formula
to derive estimates on the scalar curvature.
{\bf
\subsection{The Twisted Dirac Operator}}\label{DiracSubs}
Let~$(M,g)$ be a Riemannian manifold.
The Riemannian curvature tensor~$R^M$
induces a self-adjoint curvature operator~$\mathcal R^M$
on~$\Lambda^2TM$, such that
\begin{equation}\label{CalRDef}
	g\bigl(\mathcal R^M(e_i\wedge e_j),e_k\wedge e_l\bigr)
	=-g\bigl(R^M_{e_i,e_j}e_k,e_l\bigr)\;,
\end{equation}
where~$e_1$, \dots, $e_m$ is an orthonormal base of~$TM$.
The sign has been chosen such that all sectional curvatures of~$M$
are non-negative when~$\mathcal R^M$ is non-negative,
i.e., all eigenvalues of~$\mathcal R^M$ are~$\ge 0$.

Let~$(M,g)$ and~$(N,\bar g)$ be compact oriented Riemannian manifolds,
and let~$f\colon N\to M$ be an area-non-increasing spin map.
That is,
\begin{equation}\label{AreaNonincDef}
%
|\,v\wedge w\,|_{\bar g}\; \ge \; |\,f_*v\wedge f_*w\,|_g
\end{equation}
for all~$v$, $w\in T_qN$ and all~$q\in N$;
and the second Stiefel-Whitney classes of~$TM$ and~$TN$ are related by
\begin{equation}\label{SpinMapDef}
	w_2(TN)\;=\;f^*\bigl(w_2(TM)\bigr) .
\end{equation}

Because the total Stiefel-Whitney class is multiplicative,
condition~\eqref{SpinMapDef} is equivalent to~$w_2(TN\oplus f^*TM)=0$.
In particular,
the bundle~$TN\oplus f^*TM$ admits a spin structure.
Thus we may chose a principal bundle~$P_{\Spin_n\cdot\Spin_m}\to N$
with fiber~$\Spin_n\cdot\Spin_m=(\Spin_n\times\Spin_m)/\{\pm1\}$
that projects down to the frame bundle of~$TN\oplus f^*TM$.
Let~$\S N\otimes f^*\S M$ denote
the bundle associated to the tensor product of the spinor representations.
Note that if~$M$ is spin, so is~$N$ by~\eqref{SpinMapDef},
and we may fix compatible spin structures on~$M$ and~$N$.
Then the bundles~$\S M$ and~$\S N$ exist,
and~$\S N\otimes f^*\S M$ is precisely the bundle we have just defined.
The bundle~$\S N\otimes f^*\S M$ carries a natural Hermitian metric
and a unitary connection~$\nabla$ compatible with Clifford multiplication
by elements of~$\Cl TN\otimes f^*\Cl TM$.
We will denote Clifford multiplication with~$v\in TN$ by~$\bar c(v)$
and Clifford multiplication with~$w\in f^*TM$ by~$c(w)$.

Let~$\bar D$ be the Dirac operator on~$\S N\otimes f^*\S M\to N$,
which can locally be expressed as
	$$\bar D=\sum^n_{i=1} \, \bar c(\bar e_{i}) \, \nabla_{\bar e_i} \;,$$
in terms of an orthonormal base~$\bar e_1$, \dots, $\bar e_n$
with respect to~$\bar g$.
By the BLW formula,
\begin{equation}\label{BLWFormel}
	\bar D^2
	=\nabla^* \, \nabla+\frac{\bar\kappa}{4}
		+\frac{1}{8}\sum_{i,j=1}^n\sum_{k,l=1}^m
			g\bigl(f^*R^M_{\bar e_{i},\bar e_{j}}e_k,e_l\bigr)
			\,\bar c(\bar e_{i})\bar c(\bar e_{j})
			\otimes c(e_k)c(e_l)\;,
\end{equation}
where~$e_1$, \dots, $e_m$ is a local orthonormal base of~$f^*TM$,
$\bar\kappa$ denotes the scalar curvature of~$N$,
and~$f^*R^M$ is the curvature of the bundle~$f^*TM$.
Let us define Clifford multiplication by 2-forms by
$$
	\bar c(\bar v\wedge\bar w)=\bar c(\bar v)\,\bar c(\bar w)
		\qquad\text{and}\qquad
	c(v\wedge w)=c(v)\,c(w)\;,
$$
for~$q\in N$,
$\bar v$, $\bar w\in T_qN$
and~$v$, $w\in(f^*TM)_q$ with~$\bar g(\bar v,\bar w)=g(v,w)=0$.
If~$ \{\omega_i\} $ and~$ \{\bar\omega_j\} $ are orthonormal bases
of~$f^*\Lambda^2TM$ and~$\Lambda^2TN$,
we may rewrite equation~\eqref{BLWFormel} as
\begin{equation}\label{GuteBLWFormel}
	\bar D^2
	=\nabla^* \, \nabla+\frac{\bar\kappa}{4}
		-\frac{1}{2}\sum_{i,j}
			g\bigl(\mathcal R^M(f_*\bar\omega_j),\omega_i\bigr)
			\,\bar c(\bar\omega_j)
			\otimes c(\omega_i)\;.
\end{equation}

Let~$\kappa$ and~$\rho$ denote the scalar and Ricci curvature of~$M$.
In the rest of this section, we prove the following

\begin{Lemma}\label{Estimate}
Let~$(M,g)$ and~$(N,\bar g)$ be compact, connected, oriented
Riemannian manifolds,
and let~$f\colon M\to N$ be an area-non-increasing spin map.
Suppose that the curvature operator of~$M$ is non-negative
and that the bundle~$\S N\otimes f^*\S M$ admits a $\bar D$-harmonic spinor.
Then~$\bar\kappa\ge\kappa\circ f$ everywhere on~$N$
implies that~$\bar\kappa=\kappa\circ f$.
If moreover, $\rho>0$ and~$2\rho-\kappa<0$,
then~$\bar\kappa\ge\kappa\circ f$ implies
that~$f$ is a Riemannian submersion.
\end{Lemma}


\subsection{The Estimate}\label{EstSubs}
In this subsection,
we prove the first part of Lemma~\ref{Estimate}:
we show that~$\bar\kappa\ge\kappa\circ f$
together with the existence of a $\bar D$-harmonic spinor
implies that~$\bar\kappa=\kappa\circ f$.

We start by investigating the last term on the right hand side of
the BLW formula~\eqref{BLWFormel}.
Since we have assumed that~$\mathcal R^M$ is non-negative,
it possesses a self-adjoint square root~$L\in\operatorname{End}(\Lambda^2TM)$
such that
	$$g\bigl(\mathcal R^M\omega_i,\omega_j\bigr)
	=g\bigl(L\omega_i,L\omega_j\bigr)\;.$$
Let us write
\begin{equation}\label{LbarDef}
	{\bar L} \, \omega_k 
	:=\sum_{i}g(L\omega_k,f_*{\bar \omega}_i){\bar \omega}_i
	\quad\in\Lambda^2TN\;.
\end{equation}
Now the last term on the right hand side of~\eqref{GuteBLWFormel}
can be rewritten as
\begin{multline}\label{QuadratFormel}
	-\frac{1}{2}\sum_{i,j}
		g\bigl(\mathcal R^M(f_*\bar\omega_j),\omega_i\bigr)
		\,\bar c(\bar\omega_j)\otimes c(\omega_i)\\
\begin{split}
	&=-\frac{1}{2}\sum_{i,j,k}
		g\bigl(L(f_*\bar\omega_j),\omega_k\bigr)
		\,g\bigl(L\omega_i,\omega_k\bigr)
		\,\bar c(\bar\omega_j)\otimes c(\omega_i)\\
        &=-\frac{1}{2}\sum_k
                \bar c({\bar L}(\omega_k))\otimes c(L\omega_k)\\
	&=\frac14\sum_k\Bigl(-\bigl(\bar c(\bar L\omega_k)\otimes 1
				+1\otimes c(L\omega_k)\bigr)^2
			+\bar c(\bar L\omega_k)^2\otimes 1
			+1\otimes c(L\omega_k)^2\Bigr)
\end{split}\\
	\ge\frac14\sum_k\Bigl(\bar c(\bar L\omega_k)^2\otimes 1
			+1\otimes c(L\omega_k)^2\Bigr)\;.
\end{multline}
Here we have used
that Clifford multiplication with 2-forms is skew symmetric
and that squares of skew-symmetric endomorphisms are non-positive,
so~$-\bigl(\bar c(\bar L\omega_k)\otimes 1
+1\otimes c(L\omega_k)\bigr)^2$
is a non-negative endomorphism.

We claim that the operators~$\sum_k {\bar c}  ({\bar L} \omega_k)^2$
and~$\sum_k  c (L \omega_k)^2$ act on spinors
as multiplication by functions on~$N$,
and moreover,
\begin{equation}\label{ScalarCurvEst}
\sum_k\bar c(\bar L\omega_k)^2\ge-\frac{\kappa\circ f}2\;,
		\qquad\text{and}\qquad
	\sum_kc(L\omega_k)^2=-\frac{\kappa\circ f}2\;.
\end{equation}
The term~$-\frac12\,\kappa\circ f$ in the second statement
of~\eqref{ScalarCurvEst} arises in precisely the way
as the term~$\frac\kappa4$ in the classical BLW formula,
cf.~\cite{LM}.

The proof of the first statement is similar:
By definition of~$\bar L$ in~\eqref{LbarDef},
\begin{multline}\label{SecondScalarTerm}
	\sum_k\bar c(\bar L\omega_k)^2
	=\sum_{i,j,k}g(L\omega_k,f_*\bar\omega_i)
		\,g(L\omega_k,f_*\bar\omega_j)
		\,\bar c(\bar\omega_i)\bar c(\bar\omega_j)\\
	=\sum_{i,j}g\bigl(\mathcal R^M(f_*\bar\omega_i),f_*\bar\omega_j\bigr)
		\,\bar c(\bar\omega_i)\bar c(\bar\omega_j)\;.
\end{multline}

At this point,
we choose a local $\bar g$-orthonormal frame $\bar e_1$, \dots, $\bar e_n$
and a local $g$-orthonormal frame $e_1$, \dots, $e_m$,
such that there exists~$\mu_1$, \dots, $\mu_{\min(m,n)}\ge 0$ with
	$$f_*\bar e_i
	=\begin{cases}
		\mu_i \, e_i	&\text{if~$i\le\min(m,n)$, and}\\
		0		&\text{otherwise.}
	\end{cases}$$
This can be done by diagonalizing~$f^*g$ with respect to the metric~$\bar g$.
Then we have the orthonormal bases~$\bar e_i\wedge\bar e_j$ of~$\Lambda^2TN$
and~$e_k\wedge e_l$ of~$\Lambda^2TM$,
with
	$$f_*(\bar e_i\wedge\bar e_j)
	=\mu_i\mu_j\,e_i\wedge e_j
		\qquad\text{and~$\qquad\mu_i\mu_j\le 1$}$$
for~$1\le i<j\le\min(m,n)$,
because we have assumed~$f$ to be area-non-increasing.

We rewrite equation~\eqref{SecondScalarTerm}
in these bases,
using the definition~\eqref{CalRDef} of~$\mathcal R^M$:
\begin{multline}\label{SecondScalarEst}
	\sum_k\bar c(\bar L\omega_k)^2
	=-\sum_{i<j,k<l}
		g\bigl(R^M_{f_*\bar e_i,f_*\bar e_j}f_*\bar e_k,
				f_*\bar e_l\bigr)
		\,\bar c(\bar e_i)\bar c(\bar e_j)
		\,\bar c(\bar e_k)\bar c(\bar e_l)\\
	=-\frac14\sum_{i,j,k,l}\mu_i\mu_j\mu_k\mu_l
		\,\bigl(R^M_{ijkl}\circ f\bigr)
		\,\bar c(\bar e_i)\bar c(\bar e_j)
		\,\bar c(\bar e_k)\bar c(\bar e_l)\\
	=-\frac12\sum_{i,j}\mu_i^2\mu_j^2\,\bigl(R^M_{ijji}\circ f\bigr)
	\ge-\frac{\kappa\circ f}2\;.
\end{multline}
Here, all terms with four different indices are eliminated
by Bianchi's first identity,
while all terms with three different indices vanish for symmetry reasons.
This proves our claim~\eqref{ScalarCurvEst}.

We are now ready to prove the first statement in Lemma~\ref{Estimate}.
Assume that~$\bar\kappa\ge\kappa\circ f$.
Let~$0\ne\psi\in\Gamma(\S N\otimes f^*\S M)$ be a $\bar D$-harmonic spinor,
and let~$\left\|\mathord{\cdot}\right\|$
and~$\langle\mathord{\cdot},\mathord{\cdot}\rangle$
denote the $L^2$ norm and $L^2$ scalar product
on~$\Gamma(\S N\otimes f^*\S M)$.
Then by equations~\eqref{GuteBLWFormel}, \eqref{QuadratFormel}
and~\eqref{ScalarCurvEst},
\begin{multline}\label{FinalEstimate}
	0=\left\|\bar D\psi\right\|^2
	=\left\|\nabla\psi\right\|^2
		+\biggl\langle\psi,\biggl(\frac{\bar\kappa}4
			-\frac12\sum_{i,j}
			g\bigl(\mathcal R^M(f_*\bar\omega_j),\omega_i\bigr)
			\,\bar c(\bar\omega_j)\otimes c(\omega_i)\biggr)
			\psi\biggr\rangle\\
	\ge\left\|\nabla\psi\right\|^2
		+\biggl\langle\psi,\frac{\bar\kappa-\kappa\circ f}4
			\,\psi\biggr\rangle
	\ge 0\;.
\end{multline}
Because~$N$ is connected and~$\psi\ne0$ is $\bar D$-harmonic,
the subset of~$N$ where~$\psi$ is non-zero is open and dense in~$N$.
In particular, the estimate~\eqref{FinalEstimate}
now implies~$\bar\kappa=\kappa\circ f$.
This proves the first claim of Lemma~\ref{Estimate}.


\subsection{The Rigidity Statement}\label{RigidSubs}
We will now establish the second claim in Lemma~\ref{Estimate}.
We have to show that~${\bar \kappa} \ge \kappa\circ f$
implies that~$f$ is a Riemannian submersion
if~$\bar D$ possesses a harmonic spinor
and the Ricci curvature satisfies~$\rho>0$ and~$2\rho-\kappa<0$.

By the arguments of the last section,
${\bar \kappa} \ge \kappa\circ f$ implies
that all inequalities in~\eqref{QuadratFormel}, \eqref{SecondScalarEst}
and~\eqref{FinalEstimate} turn into equalities.
From~\eqref{SecondScalarEst}, we get in particular that
	$$\sum_{i,j}\mu_i^2\mu_j^2\,\bigl(R^M_{ijji}\circ f\bigr)
	=\kappa\circ f
	=\sum_{i,j}R^M_{ijji}\circ f\;,$$
so
\begin{equation}\label{RigidityEst}
	0=\sum_{i,j}\bigl(1-\mu_i^2\mu_j^2\bigr)
		\,\bigl(R^M_{ijji}\circ f\bigr)\;.
\end{equation}
Since $ R^M_{ijji}\ge 0 $  and~$\mu_i \mu_j\le 1$
because~$f$ is area-non-increasing,
all summands are non-negative.

Assume first that~$f$ is length-non-increasing,
i.e., $\mu_i \le 1$ for all $i$. 
Because we have assumed the Ricci curvature~$\rho$ to be positive definite,
we have~$\rho_{ii}= \sum_j R^M_{ijji}> 0$.
Thus for any~$i$ there is a~$j$ with~$ R^M_{ijji}\circ f \ne 0 $.
Hence,
$\mu_i \mu_j =1$, so~$\mu_i = \mu_j = 1$.
Since we can start with any~$i\in\{1,\dots,m\}$,
we get~$\mu_1=\dots=\mu_m=1$.
This implies in particular that~$m=\dim M\le n=\dim N$
and that~$f$ is a Riemannian submersion.

We turn to the general case,
i.e., $f$ is now only area-non-increasing.
The condition $ 2 \rho - \kappa < 0 $ implies for fixed~$k$ that
	$$2\sum_jR^M_{kjjk}<\sum_{i,j}R^M_{ijji}\;,
		\qquad\text{so}\qquad
	0 < \sum_{i\ne k,j\ne k}R^M_{ijji}\;,$$
Hence,
there is at least one pair~$(i,j)$
with~$i\ne k$, $j\ne k$, and~$R^M_{ijji} \ne 0 $.
Then we have again~$\mu_i \mu_j = 1$.
Together with~$\mu_i\mu_k\le 1$ and~$\mu_j\mu_k\le 1$,
this clearly implies~$\mu_k\le 1$.
It follows that~$f$ is length-non-increasing.
The arguments above show that~$f$ is a Riemannian submersion.
This finishes the proof of Lemma~\ref{Estimate}.
\qed

\begin{Remark}
Let~$(M,g)$ and~$(N,\bar g)$ be as in Lemma~\ref{Estimate},
and again~$\mathcal R^M\ge 0$,
but assume that~$f$ is a length-non-increasing spin map.
Suppose that~$\bar\kappa\ge\kappa\circ f$, that~$\rho>0$,
and the bundle~$\S N\otimes f^*\S M$ admits a $\bar D$-harmonic spinor.
Then~$\bar\kappa=\kappa\circ f$ everywhere on~$N$,
and~$f$ is a Riemannian submersion.
\end{Remark}

Note that $2\rho - \kappa < 0$ implies that $\dim M \ge 3$:
The condition $2\rho - \kappa < 0 $ can be rephrased
by saying that no eigenvalue of the Ricci curvature
is larger or equal than the sum of the remaining eigenvalues.
Clearly this implies the existence of at 
least three (not necessarily different) eigenvalues of $\rho$.

On the other hand, if $M$ is a locally symmetric space of compact type
and~$\dim M\ge 3$,
then the conditions~$\rho > 0 $ and $2\rho-\kappa<0$
are automatically satisfied.
Indeed,
$M$ splits locally into irreducible components of dimension~$\ge 2$ which
are Einstein. In particular, all eigenvalues of~$\rho$ are strictly positive,
and each eigenvalue has multiplicity at least~$2$.
Together with~$\dim M\ge 3$, this implies that no eigenvalue of~$\rho$
can be larger or equal than the sum of the remaining eigenvalues.

Finally, we remark that for manifolds with a non-negative curvature operator,
the two conditions~$\rho>0$ and~$2\rho-\kappa<0$ are only restrictive if
the universal cover of~$M$ contains factors which are either flat or
non-symmetric spheres or complex projective spaces.


\section{Index-Theoretic Considerations}\label{InSect}
In order to apply the results of the previous section
to a specific map~$f\colon N\to M$,
we have to ensure that the operator~$\bar D$ of Section~\ref{DiracSubs}
has a non-zero kernel.
We list some criteria that imply the existence of $\bar D$-harmonic spinors.

\subsection{Manifolds with Non-Vanishing Euler Characteristic}
In the simplest application of Lemma~\ref{Estimate},
we take a Riemannian manifold~$(M,g)$ with~$\mathcal R^M\ge 0$
and non-vanishing Euler characteristic~$\chi(M)$.
Recall that~$\mathcal M(g)$ was defined in \eqref{MgDef}.
If we take~$\bar g\in\mathcal M(g)$,
then the identity map~$\id_M\colon(M,\bar g)\to(M,g)$
is area-non-increasing and spin by~\eqref{AreaNonincDef}
and~\eqref{SpinMapDef}.

\begin{Theorem}\label{EulerSymmTh}
Let~$(M,g)$ be a compact, connected, oriented Riemannian manifold
with non-negative curvature operator 
and non-vanishing Euler characteristic.
Let~$\bar g\in\mathcal M(g)$,
then~$\bar\kappa\ge\kappa$ implies~$\bar\kappa=\kappa$.
If moreover, the Ricci curvature of~$g$ satisfies~$\rho>0$
and~$2\rho-\kappa<0$,
then~$\bar\kappa\ge\kappa$ implies~$\bar g=g$.
\end{Theorem}

\begin{proof}
Let~$\S$ denote the spinor bundle of~$M$,
which exists over all sufficiently small open subsets of~$M$
even if~$M$ is not spin.
We equip~$\S$ with the metric and connection induced by~$g$.
If the Euler characteristic~$\chi(M)$ of~$M$ is non-zero,
then~$M$ is even-dimensional,
and the local spinor bundle splits as~$\S=\S^+\oplus\S^-$.
It is a well known fact that
	$$\Lambda^\even TM
	=\bigl(\S^+\tensor\S^+\bigr)\oplus\bigl(\S^-\tensor\S^-\bigr)
		\qquad\text{and}\qquad
	\Lambda^\odd TM
	=\bigl(\S^+\tensor\S^-\bigr)\oplus\bigl(\S^-\tensor\S^+\bigr)\;,$$
that the Dirac operator on~$\Lambda^*TM=\S\tensor\S$
is precisely the operator~$D=d+d^*$,
and that the index of~$D\colon\Omega^\even(M)\to\Omega^\odd(M)$
equals~$\chi(M)$.

Let~$\bar\S$ denote the (local) spinor bundle of~$M$,
equipped with the metric and connection induced by~$\bar g$.
Then the operator~$\bar D$ considered in Section~\ref{DiracSubs}
is precisely the twisted Dirac operator on~$\bar\S\tensor\S$.
If we introduce a grading of~$\bar\S\tensor\S$
analogous to the grading of~$\Lambda^*TM$ by even and odd degree,
then the index of~$\bar D$ with respect to this grading
again equals~$\chi(M)$.
In particular,
there is a $\bar D$-harmonic spinor~$0\ne\psi\in\Gamma(\bar\S\tensor\S)$.
Now our claim follows from Lemma~\ref{Estimate}.
\end{proof}

Recall that a symmetric space~$M=G/K$ of compact type
has~$\mathcal R^M\ge 0$.
Moreover, $\chi(M)\ne 0$ iff~$\rk G=\rk K$.
Hence we have

\begin{Corollary}\label{cor}
Let $(M=G/K,\,g)$ be a compact Riemannian symmetric space  with
$\rk G=\rk K$. If~$\bar g\in\mathcal M(g)$,
then~$\bar\kappa\ge\kappa$ implies~$\bar g=g$.
\end{Corollary}

\begin{Remark}
We could also consider another grading of the bundle~$\bar\S\tensor\S$
analogous to the splitting of~$\Lambda^*TM$ into self-dual and
anti-self-dual forms.
The index of~$\bar D$ with respect to this grading
is the signature~$\operatorname{Sign}(M)$.
By Hirzebruch's signature theorem,
$\operatorname{Sign}(M)$ can be expressed
as a certain Pontrjagin number of~$M$.
A classical result of Bott implies that all Pontrjagin numbers
of a quotient~$G/K$ of compact Lie groups vanish unless~$\rk G=\rk K$.
Thus,
we do not gain anything here if we consider the signature
instead of the Euler characteristic.
\end{Remark}


\subsection{Maps of Non-Vanishing $\hat A$-Degree}
In this section,
we investigate a certain class of maps to manifolds with a non-negative
curvature operator on 2-vectors.
In order to state our result,
let us recall the following definition:
the {\em $\Adach$-degree\/} of~$f$ is given by
	$$\deg_{\Adach}f
	=\bigl(\Adach(N)\,f^*\omega\bigr)[N]\;,$$
where~$\omega\in H^m(M,\Z)$ is the fundamental class of~$M$
corresponding to the orientation of~$M$.
Recall that the notion of an area-non-increasing spin map
was defined in \eqref{AreaNonincDef} and~\eqref{SpinMapDef}.

\begin{Theorem}\label{MapTh}
Let $(M,\,g)$ be a compact connected oriented Riemannian manifold with 
non-negative curvature operator
and with non-vanishing Euler characteristic.
Let~$(N,\bar g)$ be a compact connected oriented Riemannian manifold,
and let~$f\colon N\to M$ be an area-non-increasing spin map of non-vanishing
$\Adach$-degree.
Then~$\bar\k\ge\k\circ f$ implies~$\bar\k=\k\circ f$.
If moreover, the Ricci curvature~$\rho$ of~$M$ satisfies~$\rho>0$
and~$2\rho-\kappa<0$,
then~$f\colon N\to M$ is a Riemannian submersion.
\end{Theorem}

\begin{proof}
Since~$f$ is spin,
we can construct the bundle~$\S N\tensor f^*\S M$
and the Dirac operator~$\bar D$ on~$\Gamma(\S N\tensor f^*\S M)$
as in Section~\ref{DiracSubs}.
By the Atiyah-Singer index theorem,
the index of
	$$\bar D
	\colon\bigl(\S^+N\otimes f^*S^+M\bigr)
		\oplus\bigl(\S^-N\otimes f^*S^-M\bigr)
	\longrightarrow
	\bigl(\S^-N\otimes f^*S^+M\bigr)
		\oplus\bigl(\S^+N\otimes f^*S^-M\bigr)$$
is given by
	$$\ind(\bar D)
	=\bigl(\Adach(N)\,f^*\ch(\S^+M-\S^-M)\bigr)[N]
	=\deg_{\Adach}f\cdot\chi(M)\;,$$
because~$\ch(S^+M-\S^-M)\in H^m(M,\R)$ equals the Euler class of~$M$.
Under the hypotheses of the theorem,
we have~$\ind(\bar D)\ne 0$.
Now, the theorem follows from Lemma~\ref{Estimate}.
\end{proof}


\begin{Remark}
The conditions that~$f$ be spin and~$\deg_\Adach(f)\ne0$ look very technical.
To see that they are necessary,
suppose that~$M$ is a point.
In this situation, Theorem~\ref{MapTh} becomes precisely Lichnerowicz' theorem,
which states that a compact, connected, oriented spin manifold~$N$
of non-vanishing $\hat A$-genus 
cannot carry a metric~$\bar g$ with~$\bar\kappa\ge 0$ and strict inequality
somewhere on~$N$.
This gives us a hint how to construct counterexamples to Theorem~\ref{MapTh}
without the assumptions mentioned above:
If~$N$ is the Riemannian product~$M\times\C P^k$ for~$k\ge 1$,
then clearly the projection~$f$ onto the first factor is area-non-increasing,
but~$\bar\kappa>\kappa\circ f$.
However, it is easy to see in this situation that
	$$\deg_\Adach(f)=\Adach(\C P^k)[\C P^k]$$
which vanishes if~$k$ is odd,
while the map~$f$ is spin iff~$\C P^k$ is spin,
which is not the case for even~$k$.
\end{Remark}

We do not know if the condition~$\chi(M)\ne 0$,
which we need to ensure the existence of harmonic spinors,
can be omitted entirely.
However, it can be replaced by different conditions.
Here is one possible example:

\begin{Remark}\label{SignRem}
In the proof of Theorem~\ref{MapTh},
we worked with the Dirac operator on~$N$,
twisted by the virtual bundle~$f^*(\S^+M-\S^-M)$.
We could equally well twist with only one component~$f^*\S^\pm M$,
or with the sum~$f^*\S M$.
In the latter case,
the index of the corresponding Dirac operator on~$M$ is given by
	$$\bigl(\Adach(N)\,f^*\ch(\S M)\bigr)[N]=:\deg_{\Sign}(f)\;,$$
which we will call the {\em signature degree\/} of~$f$,
because for the identity~$\id_M$,
we get the signature of~$M$:
	$$\bigl(\Adach(M)\,\ch(\S M)\bigr)[M]
	=L(M)[M]
	=\Sign(M)\;.$$
Thus,
in \ref{MapTh},
we can replace the two conditions~$\chi(M)\ne0$ and~$\deg_\Adach(f)\ne 0$
by the single condition~$\deg_{\Sign}(f)\ne 0$
to obtain another version of the theorem.
\end{Remark}


\subsection{Odd-Dimensional Manifolds}
In this section,
we present an analogue of Theorem~\ref{EulerSymmTh} for
a certain class of odd-dimensional manifolds
with non-negative curvature operator.
The idea here is to use the invariance of the mod~$2$-index
of an anti-self-adjoint real Fredholm operator
in order to find a $\bar D$-harmonic spinor.
We consider a mod~$2$-index that is related to the Kervaire semi-characteristic,
as we will explain in Remark~\ref{KervaireRem} below.

\begin{Theorem}\label{OddTh}
Let~$(M,g)$ be an odd-dimensional compact, connected, oriented
Riemannian manifold
with non-negative curvature operator on~$\Lambda^2(TM)$,
and assume that the universal covering of~$M$ is homeomorphic
to a Riemannian symmetric space~$G/K$ of compact type with~$\rk G=\rk K+1$.
Let~$\bar g\in\mathcal M(g)$,
then~$\bar\kappa\ge\kappa$ implies~$\bar\kappa=\kappa$.
Moreover, if~$\rho > 0$ and $2\rho - \kappa <0$,
then~$\bar\kappa\ge\kappa$ implies~$\bar g=g$.
\end{Theorem}

Note that our assumption on~$M$ rules out the possibility
that the universal covering of~$M$ contains a Euclidean de Rham factor.
We recall the following description of the complex spinor bundle
of not necessarily irreducible symmetric spaces:

\begin{Proposition}[\cite{Parth}, \cite{Goette}]\label{ParthProp}
Let~$M=G/K$ be a symmetric space with~$\rk G=\rk K+k$.
Then the complex spinor bundle~$\S$
is locally induced by a representation~$\sigma$
of the Lie algebra~$\mathfrak k$ of~$K$,
which splits as
	$$\sigma=2^{\textstyle\left[\frac{k}{2}\right]}_{\phantom{x}}
			\,\bigoplus_{i=1}^q\sigma_i\;,$$
where~$\sigma_1$, \dots, $\sigma_q$ are certain pairwise non-isomorphic
irreducible complex representations of~$\mathfrak k$.
\end{Proposition}

Using this proposition,
we derive a splitting of
the bundle~$\Lambda^\even TM$ of even, real exterior forms on~$M$.
We refer the reader to~\cite{BtD} and~\cite{LM}
for all technical details concerning real representations
of real semi-simple Lie algebras and real Clifford algebras.

\begin{Proposition}\label{RealProp}
Let~$M=G/K$ be a symmetric space of dimension~$m=8p+1$
such that~$\rk G=\rk K+1$.
Then the bundle~$\Lambda^\even TM$
splits as a direct sum
	$$\Lambda^\even TM=\bigoplus_{i=1}^q\mathcal E^i\;,$$
such that each~$\mathcal E^i$ is a parallel,
$G$-invariant subbundle of~$\Lambda^\even TM$
which is invariant under the natural left action of~$\Cl^\even TM$.
Moreover, for each~$\mathcal E^i$,
the space of parallel sections has real dimension~$1$,
and each parallel section is $G$-invariant.
\end{Proposition}

\begin{proof}
Recall that for~$m=8p+1$,
the real spinor representation~$\sigma_\R$ of~$\Spin_m$
acts on a real vector space~$S_\R$ of real dimension~$2^{4p}$ (\cite{LM}).
The complex spinor representation
arises as~$\sigma=\sigma_\R\tensor_\R\C$ on~$S=S_\R\tensor_\R\C$.
Complex conjugation induces a $\C$-antilinear involution~$\bar\cdot$
on~$S$, which commutes with~$\sigma$.
If we restrict~$\sigma$ and~$\sigma_\R$ to the Lie algebra~$\mathfrak k$,
then~$\bar\cdot$ descends to a $\C$-antilinear involution
on each of the irreducible subrepresentations~$\sigma_i$
of Proposition~\ref{ParthProp}.
Let~$S_{\R,i}$ be the $(+1)$-eigenspace of~$\bar\cdot$ on~$S_i$,
then~$\sigma_i=\sigma_{\R,i}\tensor_\R\C$, and
	$$\sigma_\R=\bigoplus_{i=1}^q\sigma_{\R,i}\;.$$

By~\cite{LM},
we have~$\Lambda^\even\mathfrak p\cong\Cl^\even\mathfrak p
\cong S_\R\tensor_\R S_\R$.
Setting
	$$\mathcal E^i:=G\times_K\bigl(S_\R\tensor_\R S_{\R,i}\bigr)\;,$$
we obtain the decomposition of the Proposition.
By Proposition~\ref{ParthProp} and Schur's Lemma,
the trivial $K$-isotypical component of~$\sigma\tensor\sigma_i$
is one-dimensional over~$\C$.
Arguing with complex conjugation as above,
we see that the trivial $K$-isotypical component
of~$\sigma_\R\tensor_R\sigma_{\R,i}$ is one-dimensional over~$\R$.
This implies that the space of parallel sections of~$\mathcal E^i$
is also one-dimensional and $G$-invariant.
\end{proof}

\begin{proof}[Proof of Theorem~\ref{OddTh}]
We start with the following basic case:
Assume that~$(M,g)$ is a Riemannian locally symmetric space of compact type
with~$m=\dim M=8p+1$.
Let~$e_1$, \dots, $e_m$ be a local $g$-orthonormal base of~$TM$,
and let
	$$\omega_\R=c(e_1)\dots c(e_m)\in\operatorname{End}(\Lambda^*TM)$$
denote the real Clifford volume element.
We consider the real Dirac operator
\begin{equation}\label{RealDirac}
	D_\R:=\sum_{i=1}^m\omega_\R\, c(e_i)\nabla_{e_i}\;.
\end{equation}
Note that~$\omega_\R\,c(e_i)\in\Cl^\even(TM)$ for all~$i$,
so~$D_\R$ acts on~$\Omega^\even M$.
Because the adjoint of~$\omega_\R$ equals~$\omega_\R^*=-\omega_\R$
for~$m\equiv 1$ mod~$4$,
and because~$\omega_\R$ is parallel
and commutes with Clifford multiplication~$c(e_i)$,
the operator~$D_\R$ is formally anti-self-adjoint.
Moreover, $-D_\R^2$ is equal to the Hodge-Laplacian~$(d+d^*)^2$.

Since~$M$ is of compact type,
after passing to a finite cover,
we may assume that~$M=G/K$ is simply connected.
We assume that~$\rk G=\rk K+1$.
By Proposition~\ref{RealProp},
the bundle~$\Lambda^\even M$ splits as
a direct sum of parallel sub-bundles
\begin{equation}\label{RealSplitting}
	\Lambda^\even M
	=\bigoplus_{i=1}^q \mathcal E^i
	=\bigoplus_{i=1}^q \S_\R\tensor_\R \S_{\R,i}
\end{equation}
with an obvious notation;
and for each of these sub-bundles,
the space of parallel sections has real dimension~$1$.
The operator~$D_\R$ respects this splitting.
Because~$M$ is symmetric,
a form~$\alpha\in\Omega^*M$ is $D_\R$-harmonic
iff it is $d+d^*$-harmonic iff it is parallel.
In particular,
the restriction~$D_{\R,i}$ of~$D_\R$ to~$G\times_K(S_\R\tensor_\R S_{\R,i})$
has a one-dimensional kernel.

Now let~$\bar g\in\mathcal M(g)$.
Let~$\bar\S_\R$ be the real spinor bundle of~$M$,
equipped with the metric and connection induced by~$\bar g$,
which exists over all sufficiently small open subsets of~$M$.
Then the vector bundles
	$$\bar\S_\R\tensor_\R \S_{\R,i}$$
exist globally on~$M$.
Because~$\bar\S_\R\tensor_\R\S_{\R,i}$ is a Dirac bundle
with respect to~$\bar g$,
we may define real operators~$\bar D_{\R,i}$ as in~\eqref{RealDirac}.
Then~$\bar D_{\R,i}$ is an anti-self-adjoint deformation
of the operator~$D_{\R,i}$.
For an  anti-self-adjoint real Fredholm operator,
the parity of the dimension of its kernel is invariant under deformations.
Thus~$\bar D_{\R,i}$ has an odd-dimensional
and in particular non-empty kernel.
After complexification,
the operator~$\bar D$ of Section~\ref{DiracSubs}
also has a non-empty kernel,
so the theorem follows from Lemma~\ref{Estimate} in this special case.

Next,
suppose that~$(M,g)$ is an $8p+1$-dimensional, closed, compact
Riemannian manifold with non-positive curvature operator on 2-vectors,
and that~$\bar g\in\mathcal M(g)$.
Assume that the universal covering of~$M$ is homeomorphic
to a Riemannian symmetric space~$G/K$ of compact type with~$\rk G=\rk K+1$.
We may assume that~$M$ is itself homeomorphic to~$G/K$.
Then Theorem~\ref{GMCTheorem} implies that~$M$
has the same holonomy as~$G/K$.
In particular,
we still have a splitting of~$\Lambda^*TM$ as in~\eqref{RealSplitting},
and for each of the bundles~$\mathcal E^i$,
the real dimension of the space of parallel sections is~$1$.
Because the sum of the even Betti numbers is the same for~$M$ and~$G/K$,
there are no non-parallel $D_{\R,i}$-harmonic forms
(this also follows directly because~$\mathcal R^M\ge 0$, cf.~\cite{gal}).
Now, the argument continues as above.

Finally,
assume that~$(M,g)$ is as in the theorem,
i.e., $M$ is as above, but of arbitrary odd dimension.
Then there is an even number~$r\ge 2$,
such that~$M':=M\times S^r$ is $8p+1$-dimensional for some~$p$.
Suppose that~$\bar g\in\mathcal M(g)$.
We equip~$S^r$ with its standard metric~$g_0$
and define two metrics
	$$g':=g\oplus g_0\qquad\text{and}\qquad\bar g':=\bar g\oplus g_0$$
on~$M'$.
Now, $M'$ is a Riemannian manifold of the type we have just considered.
In particular,
if we construct~$\bar D'_{\R,i}$ using the real spinor bundles of~$M'$
with respect to the metrics~$g'$ and~$\bar g'$,
then~$\bar D'_{\R,i}$ has a non-empty kernel,
so the complex Dirac operator~$\bar D'$ also has a non-empty kernel.

Nevertheless,
we cannot apply Lemma~\ref{Estimate} directly,
because in general,
$\bar g'\ge g'$ does not hold on 2-vectors, cf.~\cite{Wolfram}.
By diagonalizing~$\bar g$ with respect to~$g$,
we construct a $g$-orthonormal frame~$e_1$, \dots, $e_m$ at~$p\in M$
as in Section~\ref{EstSubs},
such that the vectors~$\bar e_1=\mu_1e_1$, \dots, $\bar e_m=\mu_me_m$
form an orthonormal base with respect to~$\bar g$
for scalars~$\mu_1$, \dots, $\mu_m\ge 0$.
We also choose a $g_0$-orthonormal frame~$e_{m+1}$, \dots, $e_{m+r}$
on~$S^r$,
and define~$\bar e_{m+1}=e_{m+1}$, \dots, $\bar e_{m+r}=e_{m+r}$ as above
with~$\mu_{m+1}=\dots=\mu_{m+r}=1$.
We note that~$R^{M'}_{ijkl}=0$ unless~$i$, $j$, $k$, $l\le m$
or~$i$, $j$, $k$, $l>m$.
Because we know that~$\mu_i\mu_j\le 1$
if~$i$, $j\le m$ or~$i$, $j>m$,
the inequality~\ref{SecondScalarEst} still holds,
so our arguments of Section~\ref{EstSubs} still show
that~$\bar\kappa'\ge\kappa'$ implies~$\bar\kappa'=\kappa'$.
Because~$\kappa'=\kappa+\kappa_0$ and~$\bar\kappa'=\bar\kappa+\kappa_0$,
where~$\kappa_0$ denotes the (constant) scalar curvature of~$S^r$,
we have proved
that~$\bar\kappa\ge\kappa$ implies~$\bar\kappa=\kappa$.

Suppose that~$\bar\kappa\ge\kappa$ and that~$\rho>0$ and~$2\rho-\kappa$.
Then as before,
the inequality~\eqref{SecondScalarEst} becomes an equality,
and the analogue of~\eqref{RigidityEst} holds for~$M'$:
\begin{equation}\label{OddRigidity}
	0=\sum_{i,j=1}^{m+r}\bigl(1-\mu_i^2\mu_j^2\bigr)\,R^{M'}_{ijji}\;.
\end{equation}
Because~$\mu_i=1$ for~$i>m$ and~$R^{M'}_{ijji}=0$ for~$i\le m<j$ or vice versa,
we only have to sum over~$1\le i$, $j\le m$,
so~\eqref{OddRigidity} turns into~\eqref{RigidityEst}.
Then the reasoning of Section~\ref{RigidSubs} shows
that the metrics~$\bar g$ and~$g$ are equal.
This finishes the proof of Theorem~\ref{OddTh}.
\end{proof}

\begin{Remark}\label{KervaireRem}
Recall that the {\em Kervaire semi-characteristic\/}~$\sigma(M)\in\Z_2$
of~$M^{4k+1}$ is defined as~$\sum_{i=0}^{2k}b_{2i}\mod 2$,
where~$b_j:=\dim H^j(M;\R)$ denotes the $j$-th Betti number over~$\R$.
If~$M^{4k+1}=G/K$ is a Riemannian symmetric space of compact type,
then~$\sigma(M)\ne 0$ iff~$\rk G=\rk K+1$ and the number~$q$
in Propositions~\ref{ParthProp} and~\ref{RealProp} is odd.
For a compact, oriented manifold of dimension~$4k+1$ with~$\sigma(M)\ne 0$,
we see immediately that the operator~$\bar D_\R$ constructed above
has a non-trivial kernel,
cf.~\cite{LM}, Example~II.7.7.
We can thus reformulate Theorem~\ref{OddTh}
for compact, connected, oriented Riemannian manifolds~$M$
with~$\sigma(M)\ne 0$ and~$\mathcal R^M\ge 0$.
However,
it is not clear if such a reformulation will give us any new example
of a compact, connected, oriented Riemannian manifold
with~$\mathcal R^M\ge 0$ that carries an area-extremal metric.
Such a new example could be of the type~$M/\Gamma$,
where~$M=G/K$ is a symmetric space of compact type with~$\rk G>\rk K+1$,
and~$\Gamma$ acts nontrivially on~$H^*(M,\R)$
(in particular, $\Gamma$ does not act as a subgroup of~$G$).
\end{Remark}

\begin{Remark}\label{OddMapRem}
One could generalize Theorem~\ref{MapTh}
to area-non-increasing spin maps~$f\colon N\to M$
between compact, connected, oriented Riemannian manifolds~$(M,g)$
and~$(N,\bar g)$ with~$\mathcal R^M\ge 0$,
where~$\dim M$ and~$\dim N$ are not necessarily even.
For certain pairs of dimensions,
the Atiyah-Singer index theorem provides a $K$-theoretic condition on~$f$
that is sufficient for the existence of a $\bar D$-harmonic spinor
in~$\Gamma(\S N\tensor f^*\S M)$,
so that one can apply Lemma~\ref{Estimate}.
In general,
this $K$-theoretic condition will not admit a reformulation
in terms of characteristic classes
(even the Kervaire semi-characteristic mentioned above
cannot be expressed in terms of characteristic classes).
\end{Remark}


\begin{thebibliography}{BtD}


\bibitem[1]{BtD}
  T.~Br\"ocker and T.~tom Dieck,
  Representations of Compact Lie Groups
  (Springer, Berlin, 1985).

\bibitem[2]{CC}
  H.D.~Cao and B.~Chow, 
  Compact K\"ahler manifolds with non-negative curvature operator,
  Invent. Math. 83 (1986), no. 3, 553--556. 

\bibitem[3]{gal}
  S.~Gallot and D.~Meyer, 
  Op\'{e}rateur de courbure et laplacien des formes
  diff\'{e}rentielles d'une vari\'{e}t\'{e} riemannienne,
  J.\ Math.\ Pures.\ Appl.~54\ (1975), 259--284.

\bibitem[4]{Goette}
  S.~Goette,
  \"Aquivariante $\eta$-Invarianten homogener R\"aume,
  (Shaker, Aachen, 1997).

\bibitem[5]{GS}
  S.~Goette and U.~Semmelmann, 
  $\Spin^c$ Structures and Scalar Curvature Estimates,
  Preprint (1999).

\bibitem[6]{gromov}
  M.~Gromov, 
  Positive curvature, macroscopic dimension, spectral gaps
  and higher signatures,
  in: S.~Gindikin, J.~Lepowski and R.~L.~Wilson.eds.,
  Functional Analysis on the Eve of the 21st Century, Vol.~II,
  Progress in Mathematics Vol.~132 (1996), 1--213.

\bibitem[7]{KW}
  J.~L.~Kazdan and F.~W.~Warner,
  Scalar curvature and conformal deformation of Riemannian structure,
  J.\ Diff.\ Geom.~10\ (1975), 113--134.

\bibitem[8]{Wolfram}
  W.~Kramer, 
  Der Dirac-Operator auf Faserungen,
  Dissertation, Bonner Mathematische Schriften~317,
  Universit\"at Bonn, 1999.

\bibitem[9]{LM}
  H.~B.~Lawson, Jr. and M.-L.~Michelsohn,
  Spin Geometry,
  (Princeton Univ.\ Press, Princeton, N.\ J.,\ 1989).

\bibitem[10]{ll1}
  M.~Llarull, 
  Scalar curvature estimates for $(n+4k)$-dimensional manifolds,
  Diff.\ Geom.\ Appl.~6\ (1996), 321--326.

\bibitem[11]{ll2}
  M.~Llarull, 
  Sharp Estimates and the Dirac Operator,
  Math.\ Ann.~310\ (1998), 55--71.

\bibitem[12]{Lohkamp}
  J.~Lohkamp,
  Scalar Curvature and Hammocks,
  Math.\ Ann.~313\ (1999), 385--407.

\bibitem[13]{MinOo}
  M.~Min-Oo, 
  Scalar Curvature Rigidity of Certain Symmetric Spaces,
  in: Geometry, Topology and Dynamics, Montreal, PQ, 1995,
  CRM Proc.\ Lecture Notes, 15,
  (Amer.\ Math.\ Soc., Providence, RI, 1998) 127--136.

\bibitem[14]{Parth}
  K.~R.~Parthasarathy, 
  Dirac operator and the discrete series,
  Ann.\ of Math.~96\ (1972), 1--30.

\bibitem[15]{T}
  S.~Tachibana,
  A theorem of Riemannian manifolds of positive curvature operator,
   Proc. Japan Acad. 50 (1974), 301--302. 

\end{thebibliography}
\end{document}